\documentclass{amsart}

\usepackage{mathtools}
\usepackage{amsmath,amssymb,amsthm,amsfonts,mathrsfs}
\usepackage[dvipsnames]{xcolor}
\usepackage{graphpap,paralist}
\usepackage[mathscr]{eucal}
\usepackage{mathabx}
\usepackage[pdftex,colorlinks,citecolor=blue]{hyperref}
\usepackage{yfonts}
\usepackage{enumerate}
\usepackage[alphabetic]{amsrefs}
\usepackage{array}

\makeatletter
\newcommand{\thickhline}{%
    \noalign {\ifnum 0=`}\fi \hrule height 1.25pt
    \futurelet \reserved@a \@xhline
}
\newcolumntype{"}{@{\hskip\tabcolsep\vrule width 1.25pt\hskip\tabcolsep}}
\makeatother

\usepackage{caption}
\usepackage{subcaption}

\usepackage{thm-restate}

\usepackage[normalem]{ulem}

\usepackage[margin=1in]{geometry}

\numberwithin{equation}{section}

\theoremstyle{definition}
\newtheorem{theorem}{Theorem}[section]

\newtheorem{conjecture}[theorem]{Conjecture}

\newtheorem*{theorem*}{Theorem}
\theoremstyle{remark}
\newtheorem{remark}[theorem]{Remark}
\newtheorem{example}[theorem]{Example}

\DeclareMathOperator{\wt}{wt}

\title[Chromatic symmetric functions of claw-free graphs are not Schur positive]{Chromatic symmetric functions of claw-free graphs are not Schur positive}

\author{Jacob P. Matherne}
\address{Department of Mathematics, North Carolina State University, Raleigh, NC., USA}
\email{jpmather@ncsu.edu}

\author{Alejandro H. Morales}
\address{LACIM, D\'epartement de Math\'ematiques, Universit\'e du Qu\'ebec \`a Montr\'eal, Montr\'eal, QC, Canada}
\email{morales\_borrero.alejandro@uqam.ca}

\begin{document}

\begin{abstract}
Chromatic symmetric functions are well-studied symmetric functions in algebraic combinatorics that generalize 
chromatic polynomials of graphs. In 1995, Stanley introduced these symmetric functions and conjectured that they are Schur positive for  claw-free graphs. We give examples of two line graphs, which are thus claw-free, whose chromatic symmetric functions have a negative coefficient in their Schur expansion. We also give a counterexample to the 2018 conjecture of Monical that Schur positive chromatic symmetric functions have saturated Newton polytope when expanded in any finite number of variables.  Both of these examples were found using ChatGPT-5.6 Sol Pro.
\end{abstract}

\maketitle

\section{Introduction} \label{sec:intro}

Given a simple graph $G$ with vertices $V(G)=\{1,2,\ldots,n\}$, the {\em chromatic symmetric function}  of $G$, denoted by $X_G({\bf x})$, is a symmetric function introduced by Stanley in \cite{StChrom1} that generalizes the chromatic polynomial of a graph.  The chromatic symmetric function of $G$ is defined as follows:
\[
X_G({\bf x}) = 
\sum_{\kappa} x_{\kappa(1)}x_{\kappa(2)} \cdots x_{\kappa(n)},
\]
where the sum is over proper colorings $\kappa\colon V(G)\to \mathbb{N}$ of $G$.

There are important open questions about the chromatic symmetric function $X_G({\bf x})$ of a graph $G$. For
instance, Stanley conjectured in \cite[Conjecture 1.4]{StChrom2} that $X_G({\bf x})$
expands positively in the Schur function basis (i.e. is $s$-positive) whenever $G$ is claw-free.  Recall that a graph is claw-free if it has no induced subgraph isomorphic to the claw graph $K_{1,3}$.

\begin{conjecture}[Stanley \cite{StChrom1}]\label{conj:stanley claw-free}
If $G$ is claw-free, then $X_G({\bf x})$ is $s$-positive.
\end{conjecture}

Stanley checked that the claw $K_{1,3}$ is not $s$-positive and verified this conjecture for {\em co-bipartite graphs}, which are complements of bipartite graphs. Gasharov proved it for claw-free incomparability graphs of posets \cite{G2}. In fact, the symmetric functions in this case expand positively in the basis of elementary symmetric functions (i.e. are $e$-positive), a stronger fact conjectured by Stanley and Stembridge in 1993 \cite{StSt}, reduced to the chordal case by  Guay-Paquet in 2013 \cite{MGP}, and verified by Hikita in 2024 via a probabilistic argument \cite{hikita2024proofstanleystembridgeconjecture}. See also another later proof by Griffin--Mellit--Romero--Weigl--Wen via operators on Macdonald polynomials \cite{GMRWW}. For the case of claw-free and chordal incomparability graph, there are other proofs of $s$-positivity using representation theory of the symmetric group \cite{BrCh,MGP2,PrecupSommers} and traces of Hecke algebras \cite{haiman,CHSS}. In contrast, Dahlberg--Foley--van Willigenburg found families of graphs avoiding claws in a different sense with non $e$-positive chromatic symmetric functions \cite{DFvW}. Additional evidence Conjecture~\ref{conj:stanley claw-free} is the proof of real-rootedness of the independence polynomial of claw-free graphs by Chudnovsky and Seymour \cite{CS07}, which Stanley showed would be a corollary if the conjecture holds. On the other hand, $s$-positivity of symmetric functions is in general a rare phenomenon \cite{Patrias}.

A related conjecture is one due to Monical \cite[Conjecture 2.37]{Mphd}, which we now describe. The Newton polytope of a  polynomial $p({\bf x}) \in \mathbb{R}[x_1,\ldots,x_k]$ is the convex hull in $\mathbb{R}^k$ of the support of $p$. The polynomial $p$ is said to be SNP (or have saturated Newton polytope) if its Newton polytope is saturated, i.e. if the support of $p$ is equal to the set of lattice points in the Newton polytope of $p$ \cite{MTY}. The study of the SNP property of chromatic symmetric polynomials $X_G(x_1,\ldots,x_k)$ started with work of Monical \cite{Mphd} and Adve--Robichaux--Yong \cite{ARY,ARY2}. 

\begin{conjecture}[Monical \cite{Mphd}]\label{conj:monical M-convex}
If $X_G$ is $s$-positive, then $X_G(x_1,\ldots,x_k)$ is SNP for any $k$.
\end{conjecture}

Already for the claw $K_{1,3}$, we have that $X_{K_{1,3}}(x_1,x_2)=x_1^3x_2+x_2^3x_1$ is not SNP (since $(1,3)$ and $(3,1)$ are in the support but not the average $(2,2)$); however, $X_{K_{1,3}}$ is not $s$-positive. In \cite[Theorem 5.8]{MMS}, the authors, together with Selover, verified Monical's conjecture for claw-free incomparability graphs and showed  in \cite[Example 7.5]{MMS} that a stronger notion than the SNP property called {\em M-convexity} does not hold for $s$-positive chromatic symmetric functions. We also showed that in order to test Monical's conjecture one needs to consider graphs with $n\geq 12$ vertices (see \cite[Section 7.4]{MMS}).  This makes it hard to find a counterexample for Conjecture~\ref{conj:monical M-convex}.

Adve--Robichaux--Yong used computational complexity arguments of coloring of claw-free graphs in \cite[Section 1.2]{ARY} to reveal a tension between Conjectures~\ref{conj:stanley claw-free} and \ref{conj:monical M-convex}. Namely, it is unlikely that both conjectures hold. In this note, we show that there is actually no such tension since, as stated, both conjectures are false. In Section~\ref{sec:schurposcounterexample}, we give  examples of two line graphs (all line graphs are claw-free) with $12$ vertices and $21$ and $22$ edges respectively (see Figure~\ref{fig:ex1}) whose chromatic symmetric functions are not $s$-positive. In Section~\ref{sec:SNPcounterexample}, we exhibit a bipartite graph with $12$ vertices and $19$ edges (see Figure~\ref{fig:ex2}) with $s$-positive chromatic symmetric function and unsaturated Newton polytope. All three of these graphs were found using ChatGPT-5.6 Sol Pro.

\subsection{Related work}  \label{sec:related work}
The counterexamples (Example~\ref{counterex1} and Example~\ref{counterex2}) to Conjecture~\ref{conj:stanley claw-free} were discovered independently at approximately the same time by Jitendra Prajapati in \cite{jitendra,Github,MO}.


\section{Counterexamples to Schur positivity}\label{sec:schurposcounterexample}

Given a graph $H=(V(H),E(H))$, the {\em line graph} of $H$, which we denote by $L(H)$, is the graph whose vertices are the edges $E(H)$ of $H$, with two edges in $E(H)$ being adjacent if they share a vertex in $V(H)$. A well-known property of line graphs is that they are claw-free \cite[Theorem 7.1.18]{West}. We use the notation of symmetric functions in \cite[Chapter 7]{EC2}. In particular, for a partition $\lambda$, we denote by $m_{\lambda}$ and $s_{\lambda}$  the monomial and Schur symmetric functions which each form a basis of the ring of symmetric functions.

\begin{figure}
\centering
    \includegraphics[scale=0.8]{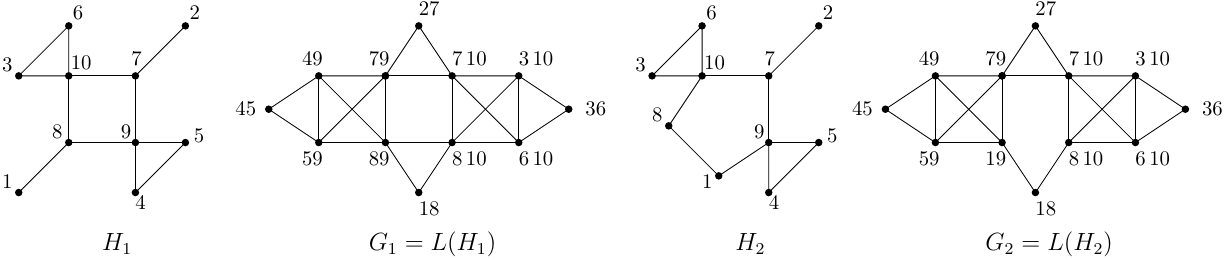}
    \caption{Graphs $H_i$ and their line graphs $G_i=L(H_i)$ for $i=1,2$. Both $G_1$ and $G_2$ are claw-free and have  chromatic symmetric functions that are not Schur positive.}
    \label{fig:ex1}
\end{figure}

\begin{example}[First counterexample to Conjecture~\ref{conj:stanley claw-free}] \label{counterex1}
Consider the graph $H_1$ with vertices in $\{1,2,\ldots,10\}$ and with $12$ edges: 
\[(1, 8), (2, 7), (3, 6), (3, 10), (4, 5), (4, 9), (5, 9), (6, 10), (7, 9), (7, 10), (8, 9), (8, 10).
\]
Let $G_1=L(H_1)$ be the line graph of $H_1$ which has $12$ vertices and $22$ edges. See Figure~\ref{fig:ex1}:Left. Since $G_1$ is a line graph, it is claw-free. Using {\tt SageMath} (see \href{https://sagecell.sagemath.org/?z=eJxVj0GLwjAQhe-C_2GOEwhhtdq6B0-CPUuPIiGEdDvQpiUTFf_9Nt2sIDl8efNmeDPNa4AjNK9hcDGQPd-9jTR6xstFrFf85ynGWdSzqIOZOrzil4RKSMCNhDJxK2Gf-Z1YSNhlHhJ3mfvFX6-wzIUyD1RZV0nf5rQpkI9Yq5680z9LrADysHy17Q2zY3XqzfMcnJsH2rTeR7uyXRgHE8lq_j9Qt_lCfGe0iliz7e5BTyNTpIdD8XYZ20-h2JreBORrIZd3E-IXcMdblQ==&lang=sage&interacts=eJyLjgUAARUAuQ==}{code}), we have that $X_{G_1}$ is not Schur positive, since
\[
[s_{3333}] X_{G_1} = -64.
\]
\end{example}

\begin{example}[Second counterexample to Conjecture~\ref{conj:stanley claw-free}]\label{counterex2}
Consider the graph $H_2$ with vertices in $\{1,2,\ldots,10\}$ and with $12$ edges: 
\[
(1, 8), (1,9), (2, 7), (3, 6), (3, 10), (4, 5), (4, 9), (5, 9), (6, 10), (7,9), (7, 10), (8, 10).
\]
Let $G_2=L(H_2)$ be the line graph of $H_2$ which has $12$ vertices and $21$ edges. See Figure~\ref{fig:ex1}:Right. Since $G_2$ is a line graph, it is claw-free. Using {\tt SageMath} (see \href{https://sagecell.sagemath.org/?z=eJxVj8EKwjAMhu-C75BjCqWoU6cHT4I7y44ipZTOBbZuNFXx7d1mFaSHrx9J-JPy1cIBylfbuhjInu7eRuo84_ks5jP-1BTjIMUgRTB9jRdcSMiFBFxKWCfuR64kbBMnzyRsEncj14mbxG3qy5Pno1-HsD6Qj1iohrzTtylVAHmYvto2htmxOjbmeQrODQPVuN1fu7J16FoTyWr-3qerdCD-MipFrNnW96D7jinSw6H4VRmrf1FsTWMC8iWT07sK8QZYrluZ&lang=sage&interacts=eJyLjgUAARUAuQ==}{code}), we have that $X_{G_2}$ is not Schur positive, since
\[
[s_{3333}] X_{G_2} = -40.
\]
\end{example}

\begin{remark}
Since Gasharov proved in \cite{G2} that the chromatic symmetric functions of claw-free incomparability graphs are $s$-positive, it follows that $G_1$ is not an incomparability graph of a poset. Indeed, one can check this since $G_1$ has a {\em net} (three vertices adjacent to a triangle) as a subgraph (e.g., $49,79,89,45,27,18$), an obstruction for being a claw-free incomparability graph  (see \cite[Figure 1]{FHM}).  A similar argument shows that $G_2$ is also not an incomparability graph of a poset.  
\end{remark}

\begin{remark}
    Combining calculations of Prajapati in \cite{jitendra} and Wang--Zhang--Zhao in \cite{WangZhangZhao} show that to isomorphism, example  $G_2$ is the smallest counterexamples to Conjecture~\ref{conj:stanley claw-free} with respect to the {\em edge-first, vertex-second} order.  
It would be interesting to find infinite families of connected claw-free graphs with chromatic symmetric functions that are not Schur positive.\footnote{After this preprint was posted, Wang, Zhang, and Zhao found such an infinite family \cite{WangZhangZhao}.}  It would also be interesting to characterize the graphs $G$ such that $X_G({\bf x})$ is $s$-positive.
\end{remark}

\section{Counterexample to the SNP property}\label{sec:SNPcounterexample}

For a coloring  \(\kappa\colon V(G) \to [k]\), define the \emph{weight} of \(\kappa\) to be
\[
\wt(\kappa) = (|\kappa^{-1}(1)|, |\kappa^{-1}(2)|, \ldots, |\kappa^{-1}(k)|) \in \mathbb{N}^k.
\]
Thus, the support of \(X_G(x_1, \ldots, x_k)\) is the set
\(
\{ \wt(\kappa) \mid \kappa \colon V(G) \to [k] \text{ is proper}\}.
\)
We call the vertex sets $\kappa^{-1}(i)$ \emph{color classes}.

\begin{figure}
\centering
    \includegraphics{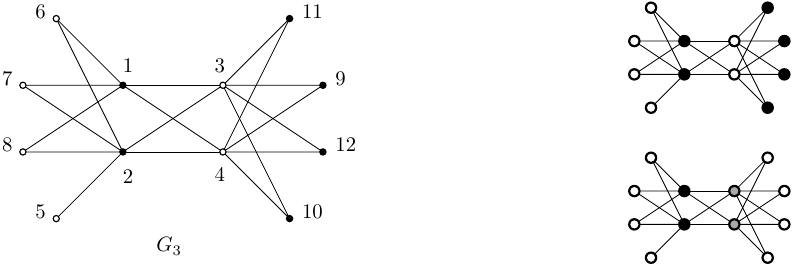}
    \caption{Bipartite graph $G_3$ that has a Schur positive chromatic symmetric function that has a non-saturated Newton polytope, since it has proper colorings of weights $(6,6,0)$ and $(8,2,2)$ (right) but not of the average $(7,4,1)$ of these weights.}
    \label{fig:ex2}
\end{figure}

\begin{example}[Counterexample to Conjecture~\ref{conj:monical M-convex}]
Let $G_3$ be the bipartite graph with vertices $\{1,2,9,10,11,12\}\cup \{3,4,5,6,7,8\}$ and with $19$ edges:
\begin{align*}
&(1,3),(1,4),(1,6),(1,7),(1,8)\\
&(2,3),(2,4),(2,5),(2,6),(2,7),(2,8)\\
&(9,3), (9,4), \quad (10,3), (10,4), \quad 
(11,3), (11,4) \quad (12,3), (12,4).
\end{align*}
See Figure~\ref{fig:ex2}. 
There are two possible proper colorings with weight $(6,6)$, namely from the bipartition of the graph: $\kappa^{-1}(i)=\{1,2,9,10,11,12\}$ and $\kappa^{-1}(j)=\{3,4,5,6,7,8\}$ for $i,j=1,2$ and $i,j=2,1$, respectively.

There are also two possible proper colorings with weight $(8,2,2)$: using one color for the $8$ vertices of degree one and two, and a different color for both vertices $\{1,2\}$, and a third color for both vertices $\{3,4\}$. See Figure~\ref{fig:ex2}:Right.

However, we claim that there is no proper coloring of weight $(7,4,1)=\frac12( (6,6,0)+(8,2,2))$, the average of the two previous weights. This is because, if there is a proper coloring $\kappa$ with such a weight then vertices $1,2,3,4$ cannot belong to the color class $\kappa^{-1}(1)$ of size $7$ since  vertices $2$, $3$, and $4$ are each not adjacent to four vertices, and vertex $1$ is not adjacent to six vertices, of which vertices $2$ and $5$ are adjacent.  Thus $\kappa^{-1}(1)$ consist of seven of the eight vertices $5,6,7,8,9,10,11,12$ (the outer vertices in the drawing in Figure~\ref{fig:ex2}). Since the vertices $1,2,3,4$ give a copy of $K_{2,2}$, at most two of these vertices have the same color. Thus these vertices and the remaining outer vertex cannot have a weight $(4,1)$, and so there is no such coloring $\kappa$.

In the monomial basis, we have shown that 
\(
[m_{660}] X_{G_3}(x_1,x_2,x_3) =  [m_{822}] X_{G_3}(x_1,x_2,x_3) = 2 \neq 0,
\)
but for the average weight $\dfrac12\left( (6,6,0)+(8,2,2)\right)=(7,4,1)$ we have that 
\(
[m_{741}] X_{G_3}(x_1,x_2,x_3)=0.
\)
This shows that $X_{G_3}(x_1,x_2,x_3)$ is not SNP. 

However, using {\tt SageMath} (see \href{https://sagecell.sagemath.org/?z=eJwLrsxVsFUIrszNTS0pykx2K81LLsnMzyvWCAzU5OWCyunlagA5HkCOe1FiQYZGtaFVtLGOiY6ZjrmORayOEYRnCuNbgvmxOoYGMIYhjAFRGlsLNC4NaJyHXnJGUX5uYklmcnwxzA3xaVBHgCwtKMrMK9FI08ssji9Ozigtii_IL84sySxL1dAEuU8jTRMAEn05Lw==&lang=sage&interacts=eJyLjgUAARUAuQ==}{code})  we have that $X_{G_3}({\bf x})$ is Schur positive since
\begin{multline*}
X_{G_2}({\bf x}) =  92528 s_{1^{12}} + 183128 s_{2 1^{10}} + 141902 s_{2^{2} 1^{8}} + 63544 s_{2^{3} 1^{6}} + 21768 s_{2^{4} 1^{4}} 
+  5348 s_{2^{5} 1^{2}} + 2612 s_{2^{6}} + 154030 s_{3 1^{9}} \\+124822 s_{3 2 1^{7}} + 52422 s_{3 2^{2} 1^{5}} 
+ 15824 s_{3 2^{3} 1^{3}} + 3319 s_{3 2^{4} 1} + 29608 s_{3^{2} 1^{6}} + 16104 s_{3^{2} 2 1^{4}} + 4625 s_{3^{2} 2^{2} 1^{2}} 
+ 1311 s_{3^{2} 2^{3}} \\+ 3018 s_{3^{3} 1^{3}} + 549 s_{3^{3} 2 1} + 142 s_{3^{4}} + 71908 s_{4 1^{8}} 
+ 49378 s_{4 2 1^{6}} + 18248 s_{4 2^{2} 1^{4}} + 4207 s_{4 2^{3} 1^{2}} + 421 s_{4 2^{4}} + 11420 s_{4 3 1^{5}} 
\\+ 6278 s_{4 3 2 1^{3}} + 1436 s_{4 3 2^{2} 1} + 1101 s_{4 3^{2} 1^{2}} + 125 s_{4 3^{2} 2} + 1272 s_{4^{2} 1^{4}} 
+ 933 s_{4^{2} 2 1^{2}} + 227 s_{4^{2} 2^{2}} + 113 s_{4^{2} 3 1} + 28 s_{4^{3}} + 20288 s_{5 1^{7}} \\+ 10660 s_{5 2 1^{5}} + 3293 s_{5 2^{2} 1^{3}} + 493 s_{5 2^{3} 1} + 1979 s_{5 3 1^{4}} + 1030 s_{5 3 2 1^{2}} 
+  185 s_{5 3 2^{2}} + 189 s_{5 3^{2} 1} + 341 s_{5 4 1^{3}} + 194 s_{5 4 2 1} \\+ 29 s_{5 4 3} 
+  69 s_{5^{2} 1^{2}} + 5 s_{5^{2} 2} + 3480 s_{6 1^{6}} + 1291 s_{6 2 1^{4}} + 321 s_{6 2^{2} 1^{2}} 
+ 54 s_{6 2^{3}} + 169 s_{6 3 1^{3}} + 90 s_{6 3 2 1} + 24 s_{6 3^{2}} + 33 s_{6 4 1^{2}} \\+ 21 s_{6 4 2} + 11 s_{6 5 1} + 2 s_{6^{2}} + 336 s_{7 1^{5}} + 82 s_{7 2 1^{3}} 
+ 22 s_{7 2^{2} 1} + 6 s_{7 3 1^{2}} + 6 s_{7 3 2} + 14 s_{8 1^{4}} + 2 s_{8 2 1^{2}} + 2 s_{8 2^{2}}.
\end{multline*}  
This shows that the graph $G_3$ is a counterexample to Monical's conjecture. 
\end{example}

\begin{remark}
Note  that the graph $G_3$ in the example above contains claws (e.g., vertices $1,6,7,8$). It would be interesting to see if there is a counterexample to Monical's conjecture that is claw-free.
\end{remark}

\section{How these examples were found}

As discussed in the introduction, while writing \cite{MMS} circa 2021, the authors became interested in the perceived tension between Conjecture~\ref{conj:stanley claw-free} and Conjecture~\ref{conj:monical M-convex}. Selover used a SAT solver in some attempts to find a counterexample to Monical's conjecture on line graphs with at least $12$ vertices. Since then, the authors of this note identified these two conjectures as amenable to be tested with LLMs. In July 2026, independently and without coordinating, the authors of this note both tried ChatGPT-5.6 Sol Pro on Monical's conjecture. For JPM's session it found  counterexample $G_3$ on the first prompt after $37$ minutes (see example of the output \href{https://chatgpt.com/s/t_6a5fc9b2c89c8191844ee7a3bc976825}{here}). For AHM's session, it did not find the counterexample after one prompt and around $98$ minutes of computation. The second prompt by AHM asked to find a counterexample to Stanley's claw-free conjecture. After $83$ minutes it found the two counterexamples $G_1$ and $G_2$ in part using SAT solvers (see example of the prompts and output \href{https://chatgpt.com/share/6a5fc49d-8fec-83ea-9a65-df6235b38c58}{here}).

\subsection*{Acknowledgements}
We thank Jesse Selover for the fruitful collaboration in \cite{MMS} that led to our study of these conjectures involving chromatic symmetric functions and for his correct intuition of using SAT solvers for them. We thank Ang\`ele Foley, Darij Grinberg, and Stephanie van Willigenburg for helpful comments, as well as Dave Anderson and  Matt Larson   for useful suggestions on how to present these counterexamples found with AI. We also thank Sam Hopkins for bringing to our attention the independent work that appeared approximately at the same time (see Section~\ref{sec:related work}).
 The second named author thanks Epoch AI for temporary access to ChatGPT Pro.  The first named author was partially supported by National Science Foundation (NSF) grant DMS-2452179 and Simons Foundation Travel Support for Mathematicians Award MPS-TSM-00007970.   The second named author was partially supported by the Natural Sciences and Engineering Research Council of Canada (NSERC) funding reference number RGPIN-2024-06246.

\bibliography{references}

@MISC {MO,
    TITLE = {Is this a counterexample to the claw-free {S}chur-positivity conjecture?},
    AUTHOR = {T. Kodama},
    HOWPUBLISHED = {MathOverflow},
    NOTE = {\url{https://mathoverflow.net/q/513515} (version: 2026-07-24)},
    EPRINT = {https://mathoverflow.net/q/513515},
    URL = {https://mathoverflow.net/q/513515}
}

@misc{Github,
  author = {J. Prajapati},
  title = {A 12-vertex counterexample to claw-free {S}chur positivity},
  year = {2026},
  publisher = {GitHub},
  journal = {GitHub repository},
  HOWPUBLISHED = {Github},
  NOTE = {\url{https://github.com/infinityscroll/claw-free-schur-counterexample}},
  URL = {\url{https://github.com/infinityscroll/claw-free-schur-counterexample}}
}

@article{hikita2024proofstanleystembridgeconjecture,
      title={A proof of the {S}tanley-{S}tembridge conjecture}, 
      author={T. Hikita},
      year={2024},
      journal={arXiv preprint arXiv:2410.12758} 
}

@ARTICLE{GMRWW,
       author = {{Griffin}, S. T. and {Mellit}, A. and {Romero}, M. and {Weigl}, K. and {Jeishing Wen}, J.},
        title = "{On Macdonald expansions of $q$-chromatic symmetric functions and the Stanley-Stembridge Conjecture}",
      journal = {arXiv e-prints},
         year = 2025,
        month = apr,
          eid = {arXiv:2504.06936},
        pages = {arXiv:2504.06936},
          doi = {10.48550/arXiv.2504.06936},
archivePrefix = {arXiv},
       eprint = {2504.06936},
 primaryClass = {math.CO},
       adsurl = {https://ui.adsabs.harvard.edu/abs/2025arXiv250406936G}
}

@article {CS07,
    AUTHOR = {Chudnovsky, M. and Seymour, P.},
     TITLE = {The roots of the independence polynomial of a clawfree graph},
   JOURNAL = {J. Combin. Theory Ser. B},
  FJOURNAL = {Journal of Combinatorial Theory. Series B},
    VOLUME = {97},
      YEAR = {2007},
    NUMBER = {3},
     PAGES = {350--357},
      ISSN = {0095-8956},
   MRCLASS = {05C69},
MRREVIEWER = {Steven D. Noble},
       DOI = {10.1016/j.jctb.2006.06.001},
       URL = {https://doi-org.libproxy.uoregon.edu/10.1016/j.jctb.2006.06.001},
}

@article {StChrom1,
    AUTHOR = {Stanley, R.~P.},
     TITLE = {A symmetric function generalization of the chromatic
              polynomial of a graph},
   JOURNAL = {Adv. Math.},
  FJOURNAL = {Advances in Mathematics},
    VOLUME = {111},
      YEAR = {1995},
    NUMBER = {1},
     PAGES = {166--194},
      ISSN = {0001-8708},
   MRCLASS = {05E05 (05C15)},
MRREVIEWER = {SeungKyung Park},
       DOI = {10.1006/aima.1995.1020},
       URL = {https://doi.org/10.1006/aima.1995.1020},
}

@incollection {StChrom2,
    AUTHOR = {Stanley, R.~P.},
     TITLE = {Graph colorings and related symmetric functions: ideas and
              applications: a description of results, interesting
              applications, \& notable open problems},
      NOTE = {Selected papers in honor of Adriano Garsia (Taormina, 1994)},
   JOURNAL = {Discrete Math.},
  FJOURNAL = {Discrete Mathematics},
    VOLUME = {193},
      YEAR = {1998},
    NUMBER = {1-3},
     PAGES = {267--286},
      ISSN = {0012-365X},
   MRCLASS = {05E05 (05C15 06A06)},
MRREVIEWER = {Daniel Ashlock},
       DOI = {10.1016/S0012-365X(98)00146-0},
       URL = {https://doi.org/10.1016/S0012-365X(98)00146-0},
}

@book {EC2,
    AUTHOR = {Stanley, R.~P.},
     TITLE = {Enumerative combinatorics. {V}olume 2},
    SERIES = {Cambridge Studies in Advanced Mathematics},
    VOLUME = {62},
 PUBLISHER = {Cambridge University Press, Cambridge},
      YEAR = {1999},
     PAGES = {xii+581},
      ISBN = {0-521-56069-1; 0-521-78987-7},
   MRCLASS = {05A15 (05-02 05E05 05E10 68R05)},
       DOI = {10.1017/CBO9780511609589},
       URL = {https://doi.org/10.1017/CBO9780511609589},
}

@article {ARY,
    AUTHOR = {Adve, A. and Robichaux, C. and Yong, A.},
     TITLE = {Computational complexity, {N}ewton polytopes, and {S}chubert
              polynomials},
   JOURNAL = {S\'{e}m. Lothar. Combin.},
  FJOURNAL = {S\'{e}minaire Lotharingien de Combinatoire},
    VOLUME = {82B},
      YEAR = {2020},
     PAGES = {Art. 52, 12},
   MRCLASS = {05E14 (05E05 14M15 14M25)},
}

@article {ARY2,
    AUTHOR = {Adve, A. and Robichaux, C. and Yong, A.},
     TITLE = {An efficient algorithm for deciding vanishing of {S}chubert
              polynomial coefficients},
   JOURNAL = {Adv. Math.},
  FJOURNAL = {Advances in Mathematics},
    VOLUME = {383},
      YEAR = {2021},
     PAGES = {Paper No. 107669, 38}
}

@article {MTY,
    AUTHOR = {Monical, C. and Tokcan, N. and Yong, A.},
     TITLE = {Newton polytopes in algebraic combinatorics},
   JOURNAL = {Selecta Math. (N.S.)},
  FJOURNAL = {Selecta Mathematica. New Series},
    VOLUME = {25},
      YEAR = {2019},
    NUMBER = {5},
     PAGES = {Paper No. 66, 37},
      ISSN = {1022-1824},
   MRCLASS = {05E05 (05E10)},
MRREVIEWER = {Allan Berele},
       DOI = {10.1007/s00029-019-0513-8},
       URL = {https://doi.org/10.1007/s00029-019-0513-8},
}

@phdthesis{Mphd,
author = {Monical, C.},
year = {2018},
title = {Polynomials in algebraic combinatorics},
school={University of Illinois Urbana-Champaign},
url={https://www.ideals.illinois.edu/bitstream/handle/2142/101462/MONICAL-DISSERTATION-2018.pdf?sequence=1&isAllowed=y}
}

@article{G2,
author = {Gasharov, V.}, 
journal = {Discrete Math.},
title = {Incomparability graphs of (3+1)-free posets are {$s$}-positive},
number = {157},
pages = {211--215},
year = {1996},
}

@article{StSt,
author = {Stanley, R. P. and  Stembridge, J. R.}, 
journal = {J. Combin. Theory Ser. A},
title = {On immanants of {J}acobi-{T}rudi matrices and permutations with restricted position},
number = {2},
pages = {261--279},
volume = {62},
year = {1993},
}

@article{MGP2,
  title={A second proof of the {S}hareshian--{W}achs conjecture, by way of a new {H}opf algebra},
  author={Guay-Paquet, M.},
journal ={arXiv preprint arXiv:1601.05498},
year={2016}
}

@article{MGP,
author={Guay-Paquet, M.},
        title = "{A modular relation for the chromatic symmetric functions of (3+1)-free posets}",
      journal = {arXiv e-prints},
         year = 2013,
        month = jun,
          eid = {arXiv:1306.2400},
        pages = {arXiv:1306.2400},
          doi = {10.48550/arXiv.1306.2400},
archivePrefix = {arXiv},
       eprint = {1306.2400},
 primaryClass = {math.CO},
       adsurl = {https://ui.adsabs.harvard.edu/abs/2013arXiv1306.2400G},
}

@article {BrCh,
    AUTHOR = {Brosnan, P. and Chow, T. Y.},
     TITLE = {Unit interval orders and the dot action on the cohomology of
              regular semisimple {H}essenberg varieties},
   JOURNAL = {Adv. Math.},
    VOLUME = {329},
      YEAR = {2018},
     PAGES = {955--1001}
}

@article{PrecupSommers,
 author = {Precup, M. and Sommers, E.},
 title = {Perverse sheaves, nilpotent {Hessenberg} varieties, and the modular law},
 fjournal = {Pure and Applied Mathematics Quarterly},
 journal = {Pure Appl. Math. Q.},
 issn = {1558-8599},
 volume = {21},
 number = {1},
 pages = {495--540},
 year = {2025},
 language = {English},
 doi = {10.4310/PAMQ.241203042708},
 zbMATH = {7964116},
 Zbl = {1565.14103}
}

@article {FHM,
    AUTHOR = {Foley, A. M. and Ho\`ang, C. T. and Merkel, O. D.},
     TITLE = {Classes of graphs with {$e$}-positive chromatic symmetric
              function},
   JOURNAL = {Electron. J. Combin.},
  FJOURNAL = {Electronic Journal of Combinatorics},
    VOLUME = {26},
      YEAR = {2019},
    NUMBER = {3},
     PAGES = {Paper No. 3.51, 19},
   MRCLASS = {05E05 (05C15)},
MRREVIEWER = {Timothy Y. Chow},
       DOI = {10.37236/8211},
       URL = {https://doi.org/10.37236/8211},
}

@article{DFvW,
 author = {Dahlberg, S. and Foley, A. and van Willigenburg, S.},
 title = {Resolving {Stanley}'s {{\(e\)}}-positivity of claw-contractible-free graphs},
 fjournal = {Journal of the European Mathematical Society (JEMS)},
 journal = {J. Eur. Math. Soc. (JEMS)},
 issn = {1435-9855},
 volume = {22},
 number = {8},
 pages = {2673--2696},
 year = {2020},
 doi = {10.4171/JEMS/974},
 zbMATH = {7227744},
 Zbl = {1444.05142}
}

@article {MMS,
    AUTHOR = {Matherne, J. P. and Morales, A. H. and Selover,
              J.},
     TITLE = {The {N}ewton polytope and {L}orentzian property of chromatic
              symmetric functions},
   JOURNAL = {Selecta Math. (N.S.)},
  FJOURNAL = {Selecta Mathematica. New Series},
    VOLUME = {30},
      YEAR = {2024},
    NUMBER = {3},
     PAGES = {Paper No. 42, 35},
      ISSN = {1022-1824,1420-9020},
   MRCLASS = {05E05 (05A20 05C15 06A07)},
MRREVIEWER = {Eric\ S.\ Egge},
       DOI = {10.1007/s00029-024-00928-4},
       URL = {https://doi.org/10.1007/s00029-024-00928-4},
}

@book{West,
 author = {West, D. B.},
 title = {Introduction to graph theory},
 isbn = {0-13-227828-6},
 year = {1996},
 publisher = {Upper Saddle River, NJ: Prentice Hall},
 language = {English},
 zbMATH = {854567},
 Zbl = {0845.05001}
}

@article{Haiman,
 author = {Haiman, M.},
 title = {Hecke algebra characters and immanant conjectures},
 fjournal = {Journal of the American Mathematical Society},
 journal = {J. Am. Math. Soc.},
 issn = {0894-0347},
 volume = {6},
 number = {3},
 pages = {569--595},
 year = {1993},
 doi = {10.2307/2152777},
 zbMATH = {279838},
 Zbl = {0817.20048}
}

@article{CHSS,
 author = {Clearman, S. and Hyatt, M. and Shelton, B. and Skandera, M.},
 title = {Evaluations of {Hecke} algebra traces at {Kazhdan}-{Lusztig} basis elements},
 fjournal = {The Electronic Journal of Combinatorics},
 journal = {Electron. J. Comb.},
 issn = {1077-8926},
 volume = {23},
 number = {2},
 pages = {56},
 note = {Id/No p2.7},
 year = {2016},
 language = {English},
 url = {www.combinatorics.org/ojs/index.php/eljc/article/view/v23i2p7},
 zbMATH = {6579082},
 Zbl = {1335.05192}
}

@article{Patrias,
 author = {Patrias, R.},
 title = {What is {Schur} positivity and how common is it?},
 fjournal = {The Mathematical Intelligencer},
 journal = {Math. Intell.},
 issn = {0343-6993},
 volume = {41},
 number = {2},
 pages = {61--64},
 year = {2019},
 language = {English},
 doi = {10.1007/s00283-018-09862-8},
 zbMATH = {7085999},
 Zbl = {1416.05299}
}

@ARTICLE{Jitendra,
       author = {{Prajapati}, J.},
        title = "{A counterexample to the claw-free Schur-positivity conjecture}",
      journal = {arXiv e-prints},
         year = 2026,
        month = jul,
          eid = {arXiv:2607.26364},
        pages = {arXiv:2607.26364},
archivePrefix = {arXiv},
       eprint = {2607.26364},
 primaryClass = {math.CO},
       adsurl = {https://ui.adsabs.harvard.edu/abs/2026arXiv260726364P}
}

@ARTICLE{WangZhangZhao,
       author = {{Wang}, D. G.~L. and {Zhang}, K. and {Zhao}, T.~Y.},
        title = "{An infinite family of counterexamples to the Stanley--Gasharov conjecture}",
      journal = {arXiv e-prints},
         year = 2026,
        month = jul,
          eid = {arXiv:2607.27166},
        pages = {arXiv:2607.27166},
archivePrefix = {arXiv},
       eprint = {2607.27166},
 primaryClass = {math.CO},
       adsurl = {https://ui.adsabs.harvard.edu/abs/2026arXiv260727166W}
}
\bibliographystyle{plain}

\end{document}